\documentstyle[11pt]{article}
\title{ Large deviations  in the reinforced random walk model on trees
} 
\author{Yu Zhang 
\\
Department of Mathematics, University of Colorado}
\begin{document}
\baselineskip .20in
\maketitle
\begin{abstract}
In this paper, we consider 
 the linearly reinforced and  the once-reinforced random walk models in the transient phase on trees. 
We  show the large deviations for the upper tails for both models.  We also show the exponential decay for the lower tail
in the once-reinforced random walk model. However,  the lower tail is in polynomial decay
for the linearly reinforced random walk model. 
\end{abstract}

\section{Introduction.}
Let ${\bf T}$ be an infinite tree with vertex set ${\bf V}$. Each $v\in {\bf V}$ has $b+1$ neighbors
except a vertex, called the {\em root}, which has $b$ neighbors for $b\geq 2$.   
We denote the root by ${\bf 0}$. For any two vertices $u,v\in {\bf V}$, let
$e=[u,v]$  be the edge with vertices $u$ and $v$. We denote by ${\bf E}$ the edge set.
Consider a Markov chain
${\bf X}=\{X_i, \omega(e, i)\}$, which starts at $X_0={\bf 0}$  with $\omega(e, 0)=1$ 
for all $e\in {\bf E}$, where $\omega(e,0)$ is called the {\em initial weight}.
For $i \geq 1$ and $e\in {\bf E}$, let $X_i\in {\bf V}$ and let $\omega(e,i)\geq 1$ be the $i$-th {\em weight}.
The transition from $X_i$ to the nearest neighbor $X_{i+1}$
is randomly selected with probabilities proportional to weights $\omega(e,i)$ of incident edges $e$ to $X_i$.
After $X_i$ has changed to $X_{i+1}$, the weights are updated by the following  rule:\\
$$
w(e, i+1)= \left\{\begin{array}{cc}
1+ k(c-1)&\mbox{ for $[X_i, X_{i+1}]=e$ and $e$ had been traversed $k$ times,}\\
 w(e, i) & \mbox{otherwise}
\end{array}
\right.
$$
for fixed $c >1$.
With this weight change, the model is called a {\em linearly reinforced random walk}.
Note that if $c=1$, then it is a simple random walk.

The linearly reinforced random walk model was first studied by Coppersmith and Diaconis in 1986 (see Diaconis (1988))
for finite graphs 
on the ${\bf Z}^d$ lattice. They asked whether the walks are recurrent or transient. For $d=1$,  the walks are recurrent for all $c\geq 1$ (see  Davis (1990) and Takeshima (2000)).   For $d\geq 1$,  Sabot and Tarres (2012) showed that the walks  are also recurrent for a large $c$. The other cases on the ${\bf Z}^d$ lattice still remain open. 
Pemantle (1988) studied this model on trees and 
showed that there exists $c_0=c_0(b)\geq 4.29$ such that when $1< c < c_0$, then the walks are transient and when $c >c_0$, 
then the walks are recurrent. Furthermore, Collevecchio (2006) and Aidekon (2008) investigated the behavior of $h(X_n)$
on the transient phase,
where $h(x)$  denotes  by the number of edges from the root to $x$ for  $x\in {\bf T}$. 
They focused on $c=2$ and showed that
the law of large numbers holds for $h(X_n)$ with a positive speed for any $b\geq 2$. 
More precisely, if $c=2$, then there exists $0< T=T(b) < b/(b+2)$ such that
$$\lim_{n\rightarrow \infty } {h(X_n)\over n} = T \mbox{ a.s.}.\eqno{(1.1)}$$
By the dominated convergence theorem,
$$\lim_{n\rightarrow \infty } {\bf E}{h(X_n)\over n} = T.\eqno{(1.2)}$$

By a simple computation, the probability that
the walks repeatedly move between an edge connected to the root  is larger than $n^{-C}$
for some $C=C(b)>0$.
Therefore,
$$n^{-C}\leq {\bf P}(h(X_n)\leq 1), \eqno{(1.3)}$$
 so the lower tail of $h(X_n)$ has the following behavior:
$$n^{-C}\leq {\bf P}(h(X_n)\leq n(T-\epsilon))\eqno{(1.4)}$$
for all $\epsilon < T$ and for all large $n$.
In this paper, $C$ and $C_i$ are positive constants depending on $c$, $b$, $\epsilon$, $N$, $M$,
and  $\delta$, but not on $n$, $m$, and $k$. They also change from appearance to appearance.
From (1.4), unlike a simple random walk on a tree, we have
$$\lim_{n\rightarrow \infty} {-1\over n^{\eta}} \log {\bf P}(h(X_n)\leq  n(T-\epsilon))=0\eqno{(1.5)}$$
for all $\epsilon < T$ and for  all $\eta >0$.

We may ask what the behavior of the upper tail is. Unlike the lower tail, we  show that the upper tail has  a standard
large deviation behavior for large $b$.\\

{\bf Theorem 1.} {\em For the linearly reinforced random walk model with $c=2$ and $b\geq 70$, and for $\epsilon >0$,
there exists a positive number $ \alpha=\alpha(b, \epsilon)$ such that
$$ \lim {-1\over n} \log {\bf P}(h(X_n)\geq  (T+\epsilon)n) =\alpha.\eqno{}$$}

{\bf Remark 1.} The proof of Theorem 1 depends on a few  Collevecchio's estimates  (see Lemma 2.1 as follows).
Since his estimates need a requirement that $b \geq 70$, Theorem 1 also needs this restriction.
We conjecture that Theorem 1 holds for all $b \geq 2$. \\

Durrett, Kesten, and Limic (2002) also investigated a similar 
reinforced random walk $\{{ Y}_k, w(e, i)\}$, except that  the weight changes by
$$
w(e, i+1)= \left\{\begin{array}{cc}
c&\mbox{ for $[Y_i, Y_{i+1}]=e$,}\\
 w(e, i) & \mbox{otherwise}
\end{array}
\right. \eqno{(1.6)}
$$
for fixed $c>1$. This random walk model is called a {\em once-reinforced random walk}.
For the once-reinforced random walk model, Durrett, Kesten, and Limic (2002) showed that for any $c>1$,
the walks are always transient.  In addition, they also showed the law of large numbers for $h(Y_n)$. 
More precisely, they showed that there exists $0< S=S(c)< b/(b+c)$ such that
$$\lim_{n\rightarrow \infty } {h(Y_n)\over n} = S\mbox{ a.s.}.\eqno{(1.7)}$$
We also investigate the large deviations for $h(Y_n)$. We have the following theorem,
similar to the linearly reinforced random walk model.\\

{\bf Theorem 2.} {\em For the once-reinforced random walk model with $c>1$ and for $\epsilon >0 $,
there exists a finite positive number $\beta=\beta(c, b, \epsilon)$ such that}
$$\lim{-1\over n} \log {\bf P}(h(Y_n)\geq  (S+\epsilon)n)=\beta. \eqno{}$$

{\bf Remark 2.} It is difficult to compute the precise  rate functions $\alpha$ and $\beta$. But we may 
obtain some properties such as the continuity in $\epsilon$ for them. \\

We may ask what the lower tail deviation for $h(Y_n)$ is. Unlike in the linearly reinforced random walk model,
the lower tail is still exponentially decaying.\\

{\bf Theorem 3.} {\em For the once-reinforced random walk model with $c>1$ and $0< \epsilon< S $,}
 $$0< \liminf {-1\over n} \log {\bf P}(h(Y_n)\leq  (S-\epsilon)n)\leq \limsup {-1\over n} 
\log {\bf P}(h(Y_n)\leq  (S-\epsilon)n)< \infty. \eqno{}$$

{\bf Remark 3.} Durrett, Kesten, and Limic (2002) also showed  that (1.7) holds  for a finitely many times reinforced random
walk.  We can also adopt the same proof of Theorems 2 and 3 to show that the same arguments hold for  a finitely many times reinforced random walk.  In fact, our proofs in Theorems 2 and  3 depend on Durrett, Kesten, and Limic's Lemmas 7 and 8  (2002). These proofs in their lemmas can be extended  for the   finitely many times reinforced random walk model.
\\

{\bf Remark 4.} We believe that the limit exists in Theorem 3, but we are unable to show it.\\


\section{Preliminaries.}  
In this section, we focus on the linearly reinforced random walk model with $c=2$.
We define a {\em hitting time} sequence $\{t_i\}$ as follows.
$${t}_k=\min\{j\geq 0: h(X_j)=k\}.$$ 
Note that walks are transient, so $h(X_j)\rightarrow \infty$ as $j\rightarrow \infty$.
Thus, $t_k$ is finite and
$$0=t_0< t_1 < t_2<\cdots < t_k<\cdots<\infty.\eqno{(2.1)}$$
With this definition, for each $k\geq 1$,
$$h(X_{{t}_k})-h(X_{{t}_{k-1}})=1.\eqno{(2.2)}$$

We also define a {\em leaving time} sequence $\{\rho_i\}$ as follows.
$$\rho_i=\max\{ j\geq 0: h(X_j)= i\}.$$
Since the walk ${\bf X}$ is transient, 
$$\rho_0 < \rho_1< \cdots < \rho_k <\cdots <  \infty.\eqno{(2.3)}$$
However, unlike the simple random walk model, $\{{t}_j-{t}_{j-1}\}$
are not independent increments. So we need to look for  independence from these times.
To achieve this target, we call $t_i$ a {\em  cut time} if
$$\rho_i-t_i=0.\eqno{(2.4)}$$

Since the walks ${\bf X}$ is transient, we may select these cut times 
and list all of them in increasing order as
$$\tau_1< \cdots < \tau_k < \cdots< \infty.\eqno{(2.5)}$$
With these cutting times, we consider diference 
$$H_k=h(X_{\tau_{k+1}})- h(X_{\tau_{k}}) \mbox{ for } k=1,2,\cdots.\eqno{(2.6)}$$
By this definition, it can be shown   that for $k=1,2,\cdots, $
$$\left(\tau_{k+1}-\tau_{k}, H_k\right) \mbox{ is an i.i.d. sequence.} \eqno{(2.7)}$$
In fact (see page 97 in Collevecchio (2006)), to verify (2.7), it is enough to realize that $X_{\tau_k}$, $k\geq 1$,
are regenerative points for the process ${\bf X}$. These points split the process ${\bf X}$ into $i.i.d.$ pieces, which
are $\{ X_m, \tau_{k} \leq m< \tau_{k+1}\}$, $k\geq 1$.

{\em Level} $k\geq 1$ is the set of vertices $v$ such that
$h(v)=k$. Level $k$ is a {\em cut} level if the walk visits it only once.
We also call $X_k$, the only  vertex to be visited, the {\em cut vertex}.
It follows from the cut time  definition that 
$X_{\tau_k}$ is a cut vertex for $k\geq 1$. 
We want to remark that 
$\tau_1$ may or may not be equal zero. If $\tau_1=0$, the root is a cut vertex. 
For convenience,  we just call  $\tau_0=0$ whether the root is a cut vertex or not. In addition, let
$$H_0=h(X_{\tau_1})-h(X_{\tau_0})=h(X_{\tau_1}).\eqno{(2.8)}$$
With these definitions, Collevecchio (2006) proved the following lemma.\\

{\bf Lemma 2.1.}  {\em For $c=2$ and $b\geq 70$, 
$${\bf P} ( H_k \geq k)\leq 0.115^k \mbox{ for } k\geq 0.\eqno{(2.9)}$$
Furthermore, for $p_0=1002/1001$,}
$${\bf E}\tau_1^{p_0} <\infty.\eqno{(2.10)}$$
With Lemma 2.1, we can see that  $h(X_{\tau_{k+1}})-h(X_{\tau_{k}})$ is large with a small probability.
Also, $\tau_{k+1}-\tau_{k}$ is large with a small probability. However, to show a large deviation result,
we need a much shorter tail requirement. Therefore, we need to truncate both 
$H_k=h(X_{\tau_{k+1}})-h(X_{\tau_{k}})$ and $\tau_{k+1}-\tau_{k}$.
We call $\tau_k$  $N$-{\em short} for $k\geq 1$ if
$$H_k=h(X_{\tau_{k+1}}) -h(X_{\tau_{k}}) \leq N;\eqno{(2.11)}$$
otherwise, we call it $N$-{\em long}.
Since we only focus on the transient phase, we have
$$\tau_k(N) < \infty. \eqno{}$$
We list all $N$-short cut times as
$$\tau_1(N)< \tau_1(N)< \cdots<\infty.\eqno{(2.12)}$$
For convenience, we also call $\tau_0(N)=0$ whether the root is a cut vertex or not.
We know that $\tau_k(N)=\tau_i$ for some $i$.
We denote it by $\tau_k'(N)=\tau_{i+1}$.  In particular, let $\tau'_0(N)=0$.
For $N>0$, let 
$$I_n=\max\{i: \tau_i(N)\leq n\}$$
 and 
$$h_n(N)=\sum_{i=0}^{I_n} \left(h\left(X_{\tau_{i}'(N)}\right)-h\left(X_{\tau_i(N)}\right)\right).\eqno{}$$
If $I_n=0$,
$$h_n=0.\eqno{(2.13)}$$

Now we state  standard tail estimates for an i.i.d. sequence. The proof can be followed directly from
Markov's inequality.\\

{\bf Lemma 2.2.} {\em Let $Z_1, \cdots Z_k,\cdots $ be an i.i.d. sequence 
with ${\bf E}Z_1=0$ and ${\bf E}\exp(\theta Z_1) < \infty$
for some $\theta >0$, and let
$$S_m=Z_1+Z_2+\cdots + Z_m.$$
For any $\epsilon >0$,  $i\leq n$ and $j\geq n$, there exist $C_i=C_i(\epsilon)$ for $i=1,2$ such that
$${\bf P}( S_i\geq n\epsilon)\leq C_1 \exp(-C_2 n),$$
and}
$${\bf P}( S_j \leq -\epsilon n)\leq C_1 \exp(-C_2 n).$$
\\

Now we show that  $h_n(N)/n$ and $h(X_n)/n$ are not very different if $N$ is large.\\

{\bf Lemma 2.3.} {\em For $\epsilon >0$, $c=2$, and  $b\geq 70$,  there exist $N=N(\epsilon)$ 
and $C_i=C_i(\epsilon,N)$ for $i=1,2$ such that}
$${\bf P}( h(X_n)\geq h_n(N) +n\epsilon)\leq C_1\exp(-C_2 n) .\eqno{}$$

{\bf Proof.}  If 
$$h(X_n)-h_n(N)\geq \epsilon n, \eqno{(2.14)}$$
we may suppose that there are only $k\geq 1$ many $N$-long  cut time pairs $\{\tau_{i_j}, \tau_{i_j+1}\}$
for $j=1,\cdots, k$ such that
$${\tau_{i_1}}<  \tau_{i_1+1}< {\tau_{i_2}}<\tau_{i_2+1}<\cdots< \tau_{{i_j}}<\tau_{i_j+1} \cdots <
\tau_{i_{k-1}}
<\tau_{i_{k-1}+1}
<\tau_{i_{k}}
\leq n\leq \tau_{{i_k}+1}$$
 with $i_1\geq 1$ and with
 $$\sum_{j=1}^{k} H_{i_j}=\sum_{j=1}^{k} h(X_{\tau_{i_{j}+1}})-h(X_{\tau_{i_j}}) \geq \epsilon n/2,\eqno{(2.15)}$$
where
$$ H_{i_j}=h(X_{\tau_{i_{j}+1}})-h(X_{\tau_{i_j}}) > N\mbox{ for } j=1,2,\cdots, k\leq n/N,\eqno{(2.16)}$$
or $$h(X_{\tau_1})\geq \epsilon n/2.\eqno{(2.17)}$$
For the second case in (2.17), by Lemma 2.1, there exist $C_i=C_i(\epsilon)$ for $i=1,2$ such that
$${\bf P} (h(X_{\tau_1})\geq \epsilon n/2)={\bf P}(H_0\geq \epsilon n/2)\leq C_1\exp(-C_2 n).\eqno{(2.18)}$$
We focus on the first case in (2.15).
By (2.7) and Lemma 2.1, $\{H_1, H_2,\cdots\}$ is an i.i.d sequence with
$${\bf P}(H_i\geq m )\leq 0.115^m \mbox{ for } i\geq 1.\eqno{(2.19)}$$
Thus, if (2.15) holds, by (2.15) and (2.16), it implies that there exist $k$ many $H_i$s in
$\{H_1,\cdots, H_{n}\}$ 
for $1\leq k\leq\lceil n/N\rceil$ such that
$H_i > N$ and their sum is large than $\epsilon n/2$.

For a fixed $k$,
it costs at most $n\choose k$ to fix the subsequence of these $H_i$s from $\{ H_1,\cdots, H_{n}\}$.
We denote by $H_{i_1}, \cdots , H_{i_k}$ these fixed random variables.  Since $\{H_i\}$
is an i.i.d sequence, the joint distribution  of $H_{i_1}, \cdots , H_{i_k}$ is always the same for different 
$i_j$s.
With these observations,
$$
{\bf P}\left(h(X_n)\geq h_n(N) +n\epsilon/2, \mbox{ (2.15) holds}\right)\leq \sum_{k=1}^{ \lceil n/N\rceil} {n\choose k} {\bf P} (H_{i_1}+\cdots +H_{i_k} \geq n\epsilon/2).
\eqno{(2.20)}$$
By (2.19), we know that
$$EH_i= EH_1< \infty \mbox{ for each }i\geq 1.$$
Since $k \leq n/N+1$, we  may take $N=N(\epsilon)$ large such that for each $k\leq n$ and fixed $i_1,\cdots, i_k$
$$ {\bf P} (H_{i_1}+\cdots +H_{i_k} \geq n\epsilon/2)\leq {\bf P} ([H_{i_1}-EH_{i_1}]+\cdots +[H_{i_k}-EH_{i_k}] \geq n\epsilon/4)
\eqno{(2.21)}$$
Note that $\{H_{i_j}-EH_{i_j}\}$ is an i.i.d sequence with a zero-mean and an exponential tail for $j=1,\cdots, k$,  so by Lemma 2.2,
$$ {\bf P} ([H_{i_1}-EH_{i_1}]+\cdots +[H_{i_k}-EH_{i_1}] \geq n\epsilon/4)\leq C_3\exp(-C_4 n).
\eqno{(2.22)}$$
By a standard entropy bound, as given in Corollary 2.6.2 of Engel (1997), for $k\leq n/N$,
$${n\choose k}\leq \exp(n\log N/N).\eqno{(2.23)}$$
By (2.19)-(2.22), if we take $N$  large, then there exist $C_i=C_i(\epsilon, N)$ for $i=5,6$ such that
$${\bf P}\left(h_n(X_n)\geq h_n(N) +n\epsilon, \mbox{ (2.15) holds}\right)\leq C_5 n\exp(-C_6 n).
\eqno{(2.24)}$$
So Lemma 2.3 holds by (2.18) and (2.24). $\Box$\\

We also need to control the time difference such that
$\tau_k'(N)-\tau_k(N)$ cannot be large. We call $\tau_k(N)$  $M$-{\em tight} for $k\geq 1$ if
$$\tau_{k}'(N) -\tau_{k}(N) \leq M.\eqno{}$$
We list all $M$-tight $N$-short cut times as 
$$\tau_1(N,M), \tau_2(N,M),\cdots, \tau_k(N,M),\cdots. \eqno{}$$
Suppose that $\tau_k(N,M) < \infty$. We know that $\tau_k(N,M)=\tau_i$ for some $i$.
We denote $\tau_k'(N,M)=\tau_{i+1}$.  
For convenience,  we also call  $\tau_0(N,M)=0$ and $ \tau_0'(N,M)=0$
whether the root is a cut vertex or not.
Let 
$$J_n=\max\{i:\tau_i(N,M)\leq n\}\eqno{}$$
and 
$$h_n(N,M)=\sum_{i=0}^{J_n} \left(h\left(X_{\tau_{i}'(N,M)}\right)-h\left(X_{\tau_i(N,M)}\right)\right).\eqno{(2.25)}$$
If $J_n=0$, then
$$h_n(N,M)=0.\eqno{(2.26)}$$

The following lemma shows that $h_n(N,M)/n$ and $h_n(N)/n$ are not far away.\\

{\bf Lemma 2.4.} {\em For $\epsilon >0$, for $N$, and for each $n$,  there exists $M=M(\epsilon, N)$ such that}
$$ h_n(N)\leq h_n(N,M) +n\epsilon .\eqno{}$$
 
{\bf Proof.} 
If $h_n(N)>h_n(N,M) +n\epsilon$, we know that there are at least $\epsilon n/2N$ many $\{\tau_i(N)\}$
such that 
$$\tau_i'(N)-\tau_{i}(N) > M.\eqno{(2.27)}$$
If we take $M \geq 3 N\epsilon ^{-1} $, then
$$n\geq \sum_{i=1}^{I_n} \left(\tau_i'(N)-\tau_{i}(N)\right) > M \epsilon n/2N > n.\eqno{(2.28)}$$
The contradiction shows that 
$$ h_n(N)\leq h_n(N,M) +n\epsilon .$$
So Lemma 2.4 follows. $\Box$\\

Let ${\cal E}(\epsilon)$ be the event that $h(X_n)\geq n(T-\epsilon)$.
By Lemmas 2.3 and 2.4,  
$$\lim_{n\rightarrow \infty}{\bf P}(h_n(N,M) \leq Tn/2, {\cal E}_n(\epsilon))=0 .$$
Note that ${\bf P}({\cal E}_n(\epsilon))$ is near  one for large $n$, so there are at least
$Tn/2M$ many $\tau_i(N,M)$s with $\tau_i(N,M) \leq n$  that also have a probability near one for large $n$.
Hence, $\tau_k(N,M)=\infty$ cannot have a positive probability for each $k$.
Therefore, 
$$\tau_1(N,M)<  \tau_2(N,M)< \cdots \tau_k(N,M)< \cdots< \infty. \eqno{(2.29)} $$
By (2.29), we know that $\tau_{k}(N,M)=\tau_i$ for some $i$ and 
$$\tau_k'(N,M)-\tau_k(N,M)=\tau_{i+1}-\tau_{i}.$$
Therefore, by the same proof of (2.7), for $k\geq 1$
$$\left\{\left(\tau_{k}'(N,M)-\tau_{k}(N,M), h\left(X_{\tau_{k}'(N,M)}\right)-h\left(X_{\tau_k(N,M)}\right)\right)\right\} \mbox{ is an i.i.d. sequence.}\eqno{(2.30)}$$

\section{Large deviations for $h_n(N,M)$.}
By Lemma 2.1, we let
$${\bf E} (\tau_2-\tau_{1})=A\geq 1\mbox{ and } {\bf E} \left(\tau_{1}'(N,M)-\tau_{1}(N,M)\right)=A(N,M) $$
and
$${\bf E}(h\left(X_{\tau_{2}})-h(X_{\tau_{1}})\right)=B\geq 1\mbox{ and }{\bf E}\left(h\left(X_{\tau_{1}'(N,M)}\right)-h\left(X_{\tau_{1}(N,M)}\right)\right)=B(N,M).$$
We set 
$$T_n=\sum_{k=1}^n (\tau_{k+1}-\tau_{k}) \mbox{ and }T_n(N,M)= \sum_{k=1}^n \left(\tau_k'(N,M)-\tau_k(N,M)\right)$$
and 
$$H_n=\sum_{k=1}^n \left(h\left(X_{\tau_{k+1}}\right)-h\left(X_{\tau_{k}}\right)\right)\mbox{ and }H_n(N,M)= \sum_{k=1}^n \left(h\left(X_{\tau_{k}'(N,M)}\right)-h\left(X_{\tau_k(N,M)}\right)\right).$$
By the law of large numbers, 
$$\lim_{n\rightarrow \infty}{T_n \over n}=A\mbox{ and }\lim_{n\rightarrow \infty}{T_n(N,M) \over n}=A(N,M)\eqno{(3.1)}$$
and 
$$ \lim_{n\rightarrow \infty}{H_n \over n}=B\mbox{ and }\lim_{n\rightarrow \infty}{H_n(N,M) \over n}=B(N,M) .\eqno{(3.2)}$$


If $\tau_i\leq n\leq \tau_{i+1}$ for $i\geq 1$, then 
$$h(X_{\tau_i})\leq h(X_n) \leq h(X_{\tau_{i+1}}).\eqno{(3.3)}$$
Thus,
$${h(X_{\tau_i})\over \tau_{i+1}}\leq {h(X_n)\over n}\leq {h(X_{\tau_{i+1}})\over \tau_i}.\eqno{(3.4)}$$
By (3.1) and (3.2),
$$\lim_{i\rightarrow \infty} {h(X_{\tau_i})\over \tau_{i+1}}=\lim_{i\rightarrow \infty} {h(X_{\tau_{i+1}})\over \tau_{i}}
={B\over A}.\eqno{(3.5)}$$
So by (1.1), (3.4), and (3.5),
$${B\over A}=T.\eqno{(3.6)}$$
Regarding $B(N,M)$ and $A(N,M)$, we have the following lemma.\\

{\bf Lemma 3.1.} {\em For $c=2$ and $b \geq 70$,}
$$\lim_{N,M\rightarrow \infty}A(N,M)=A \mbox{ and } \lim_{N,M\rightarrow \infty}B(N,M)=B
\mbox{ and } \lim_{N,M\rightarrow \infty}{B(N,M) \over A(N,M)}=T.\eqno{}$$

{\bf Proof.}  By (2.5) and the definitions of $\tau_1(N)$ and $\tau_1(N,M)$,
for each sample point $\omega$, there exist large $N$ and $M$ such that
$$\tau_1(N,M)(\omega)=\tau_1(\omega),$$
where $\tau_1(N,M)(\omega)$ and $\tau_1(\omega)$ are $\tau_1(N,M)$ and $\tau_1$
with $\omega$.
It also follows from the definition of $\tau'_1(N,M)$ that for the above $N$ and $M$,
$$\tau_1'(N,M)(\omega)=\tau_2(\omega).$$
Thus, for each $\omega$
$$\lim_{N,M \rightarrow \infty} \tau_1'(N,M)(\omega)-\tau_1(N,M)(\omega)=\tau_2(\omega)-\tau_1(\omega).\eqno{(3.7)}$$
By the dominated convergence theorem,
$$\lim_{N,M \rightarrow \infty}A(N,M)=\lim_{N,M \rightarrow \infty} {\bf E} (\tau_1'(N,M)-\tau_1(N,M))={\bf E}(\tau_2-\tau_1)=A.\eqno{(3.8)}$$
Similarly,
$$\lim_{N,M \rightarrow \infty}B(N,M)=\lim_{N,M \rightarrow \infty} {\bf E} \left(h\left(X_{\tau_1'(N,M)}\right)-h\left(X_{\tau_1(N,M)}\right)\right)=B.\eqno{(3.9)}$$
Therefore, Lemma 3.1 follows from (3.8), (3.9), and (3.6). $\Box$\\

Now we show that $h_n(N,M)$ has an exponential upper tail.\\

{\bf Lemma 3.2.} {\em If $c=2$ and $b\geq 70$, then for $\epsilon >0$, there exist $N_0=N_0(\epsilon)$ and $M_0
=M_0(\epsilon)$  such that for all $N \geq N_0$ and $M\geq M_0$} 
$${\bf P} ( h_n(N,M)\geq n (T+\epsilon))\leq C_1 \exp(-C_2 n),\eqno{(3.10)}$$
where $C_i=C_i(\epsilon, N,M)$ for $i=1,2$ are constants.\\

{\bf Proof.} Recall that
$$J_n=\max\{i:\tau_i(N,M)\leq n\}.$$
So
\begin{eqnarray*}
&&{\bf P} ( h_n(N,M)\geq n (T+B\epsilon))\\
&=& {\bf P} \left ( \sum_{i=1}^{J_n} \left(h\left(X_{\tau_{i}'(N,M)}\right)-h\left(X_{\tau_{i}(N,M)}\right)\right)\geq n(T+B\epsilon)\right)\\
&\leq & {\bf P} \left ( \sum_{i=1}^{J_n} \left(h\left(X_{\tau_{i}'(N,M)}\right)-h\left(X_{\tau_{i}(N,M)}\right)\right)\geq n(T+B\epsilon), J_n\leq 
n\left({T \over B(N,M)}+\epsilon/2\right)\right)\\
&&+ {\bf P} \left (  J_n> n\left({T \over B(N,M)}+\epsilon/2\right)\right)\\
&\leq & {\bf P} \left ( \sum_{i=1}^{n({T /B(N,M)}+\epsilon/2)} \left(h\left(X_{\tau_{i}'(N,M)}\right)-h\left(X_{\tau_{i}(N,M)}\right)\right)\geq n(T+B\epsilon)\right)\\
&&+ {\bf P} \left (  J_n> n\left({T \over B(N,M)}\right)+\epsilon/2)\right)\\
&=& I+II. \hskip 12cm (3.11)
\end{eqnarray*}
Here without loss of generality,  we assume that $n({T /B(N,M)}+\epsilon/2)$ is an integer,
otherwise we can use $\lceil n({T /B(N,M)})+\epsilon/2)\rceil$ to replace $n({T /B(N,M)})+\epsilon/2)$. 
We will estimate $I$ and $II$ separately. For $I$, note that by Lemma 3.2, there exist $N_0=N_0(\epsilon)$ and $M_0=M_0(\epsilon)$ such that
for all $N\geq N_0$ and $M\geq M_0$
$${\bf E} \left(\sum_{i=1}^{n({T /B(N,M)}+\epsilon/2)} \left(h\left(X_{\tau_{i}'(N,M)}\right)-h\left(X_{\tau_{i}(N,M)}\right)\right)\right)\leq
nT(1+B(N,M)\epsilon/2)\leq nT(1+2B\epsilon/3).$$
Note also that by (2.30),
$$\left\{h\left(X_{\tau_{i}'(N,M)}\right)-h\left(X_{\tau_{i}(N,M)}\right)\right\}\mbox{ is a uniformly bounded  i.i.d. sequence},$$
 so by Lemma 2.2,
there exist $C_i=C_i(\epsilon, N,M)$  for $i=3,4$ such that
$${\bf P} \left ( \sum_{i=1}^{n({T /B(N,M)}+\epsilon/2)} \left(h\left(X_{\tau_{i}'(N,M)}\right)-h\left(X_{\tau_{i}(N,M)}\right)\right)\geq
n(T+B\epsilon)\right)\leq 
C_3 \exp(-C_4 n).\eqno{(3.12)}$$

Now we estimate $II$. By Lemma 3.1, there exist $N_0=N_0(\epsilon, b)$ and $M_0=M_0(\epsilon, b)$ such that
for all $N\geq N_0$ and $M \geq M_0$ 
$${\bf P} \left(  J_n> n\left({T \over B(N,M)}\right)+\epsilon/2\right)={\bf P}  \left(  J_n> n\left(A^{-1}(N,M)+\epsilon/3\right)\right).\eqno{(3.13)}$$
Here without loss of generality,  we also assume that $n(A^{-1}(N,M)+\epsilon/3)$ is an integer,
otherwise we can use $\lceil n(A^{-1}(N,M)+\epsilon/3)\rceil$ to replace $n(A^{-1}(N,M)+\epsilon/3)$.
Note that
$$ \left\{J_n \geq n(A^{-1}(N,M)+\epsilon/3)\right\}\subset \left\{\sum_{i=1}^{n(A^{-1}(N,M)+\epsilon/3)} (\tau_i'(N,M)- \tau_i(N,M)) \leq n\right\}.\eqno{(3.14)}$$
Note also that
$${\bf E} \sum_{i=1}^{n(A^{-1}(N,M)+\epsilon/3)} \left(\tau_i'(N,M)- \tau_i(N,M)\right)= n(1+ \epsilon A(N,M)/3),$$
and, by (2.30), $\{\tau_i'(N,M)- \tau_i(N,M)\}$ is a uniformly bounded i.i.d. sequence, so by (3.13), and (3.14), and Lemma 2.2,
there exist $C_i=C_i(\epsilon, b,N,M)$  for $i=5,6$ such that
\begin{eqnarray*}
&& {\bf P} \left(  J_n> n\left({T \over B(N,M)}+\epsilon/2\right)\right)\\
&\leq &{\bf P}  \left(  J_n> n(A^{-1}(N,M)+\epsilon/3)\right)\\
&\leq& {\bf P}\left( \sum_{i=1}^{n(A^{-1}(N,M)+\epsilon/3)} (\tau_i'(N,M)- \tau_i(N,M)) \leq n \right)\\
&\leq & C_5\exp(-C_6 n).\hskip 11cm (3.15)
\end{eqnarray*}
For all large $N$ and $M$, we substitute (3.12) and (3.15) in (3.11) to have
$${\bf P} ( h_n(N,M)\geq n (T+\epsilon))\leq I + II \leq C_7 \exp(-C_8 n)\eqno{(3.16)}$$
for $C_i=C_i(\epsilon, N,M)$  for $i=7,8$.
Therefore, we have an exponential tail estimate  for $h_n(N,M)$. So Lemma 3.2 follows. $\Box$\\

Let
$$L_n=\max\{i: \tau_i\leq n\}$$
and 
%
$$h_n=\sum_{i=1}^{L_n}\left(h(X_{\tau_{i}})-h(X_{\tau_{i-1}})\right) \mbox{ if } L_n\geq 1 \mbox{ and } h_n=0 \mbox{ if } L_n=0.\eqno{(3.17)}$$

Recall that $\rho_i$ is the leaving time defined in (2.3).  We show the following subadditive argument for $h_n$.\\

{\bf Lemma 3.3.} {\em For $c=2$,  $b \geq 2$, $N >0$, and for each pair of positive integers $n$ and $m$, 
$${\bf P} ( h_n\geq n C, \rho_0 \leq N){\bf P} ( h_m\geq m C, \rho_0\leq N)\leq 2^N(b+1)n{\bf P} ( h_{n+m+1}\geq (n+m)  C+1,\rho_0\leq N),$$
for any $C>0$.}\\


{\bf Proof.} 
By the definition in (3.17),  there exists
$0\leq k\leq n$ such that
$$\tau_k \leq n\leq \tau_{k+1}.\eqno{}$$
So
$$h_n= h(X_{\tau_k})\leq h(X_n)\leq h(X_{\tau_{k+1}}).\eqno{(3.18)}$$
For $ i\geq nC$, we denote by ${\cal F}(x,  i, N, nC)$ the event that  walks $\{X_1, X_2, \cdots, X_{i}\}$ have 
$$h(X_j)<  nC \mbox{ for $j< i$ and } h(X_i)=x\mbox{ with } h(x) \geq nC.\eqno{(3.19)}$$
In addition, the number of   walks $\{X_1, X_2, \cdots, X_{i}\}$ visiting  the root is no more than $N$.

Note that on $\{h_n\geq n C, \rho_0\leq N\}$, walks eventually move to some vertex $x$ at some time $i$ with $h(x) \geq nC$,
 and   walks $\{X_1, X_2,\cdots, X_i\}$ visit the root  no more than $N$ times.
So we may control   $\{h_n\geq n C, \rho_0\leq N\}$ by a finite step walks $\{X_1, X_2, \cdots, X_{i}\}$
in order to work on a further coupling process. More precisely, 
$${\bf P} ( h_n\geq n C, \rho_0\leq N)\leq \sum_{i\leq n}\sum_{x} {\bf P} \left({\cal F}(x,  i, N, nC) \right).\eqno{(3.20)}$$
 There are $b+1$ many  vertices adjacent to $x$.  We just select one of them and denote it by $z$ with $h(z)=h(x)+1$. Let $e_z$ be the edge with the vertices $x$ and $z$.
  On ${\cal F}(x,i,N, nC)$, we require that the next move $X_{i+1}$ will be from $x$ to $z$. Thus,
$X_{i+1}=z$. We denote this subevent by ${\cal G}(x,z, i, N, nC)\subset {\cal F}(x,i,N, nC)$.
We have
$$\sum_{i\leq n}\sum_{x} {\bf P} \left({\cal F}(x,  i, N, nC) \right)\leq (b+1)\sum_{i\leq n}\sum_{x} {\bf P} \left({\cal G}(x,z,  i, N, nC )\right).\eqno{(3.21)}$$

Now we focus on $\{h_m\geq Cm, \rho_0\leq N\}$.
Let ${\bf T}_z$ be the subtree with the root at $z$ and vertices in $\{v: h(v) \geq h(z)\}$.
We define $\{X_n^i(z)\}$ to be the linearly reinforced random walks starting from $z$  in  subtree  ${\bf T}_z$ for $n\geq i+1$ with 
$$X_{i+1}^{i}(z)=z \mbox{ and } w(e_z, i+1)=2.$$
Note that  walks $\{X_n^i(z)\}$ stay inside ${\bf T}_z$, so
$$w(e_z, n) =2 \mbox{ for } n\geq i+1.\eqno{(3.22)}$$
We can define $\tau_k^i$, $\rho^i_0$ and $h_m^i(z)$ for  $\{X_n^i(z)\}$ similar to the definitions of $\tau_k$, $\rho_0$  and $h_m$ for $\{X_n\}$. 

 On $w(e_z,i+1)=2$,
we consider  a probability difference between  ${\bf P}(h_m \geq Cm, \rho_0\leq N)$ and ${\bf P}(h_m^i(z) \geq Cm, \rho_0^i\leq N)$.  
Note that there are only $b$ edges from the root, but there are $b+1$ edges from vertex $z$ with  $w(e_z,n)=2$, so the two probabilities are not the same.
We claim that 
$$
{\bf P}(h_m\geq m C, \rho_0\leq N)\leq 2^N {\bf P}(h_m^i(z) \geq Cm, \rho_0^i\leq N\,\,|\,\, w(e_z, i+1)=2).\eqno{(3.23)}$$
To show (3.23), we consider  a fixed path  $(u_0={\bf 0}, u_1, u_2, \cdots)$  in ${\bf T}$ with 
$\{X_1=u_1, X_2=u_2, \cdots\}\in \{h_m \geq Cm, \rho_0\leq N\} $. Note that  $[u_j, u_{j+1}]$
is an edge in ${\bf E}$.
If we remove  ${\bf T}$ from the root to $z$, it will be ${\bf T}_z$. So path $({\bf 0}, u_1, u_2, \cdots)$ in ${\bf T}$ will be a new path $(u_0(z)=z, u_1(z), u_2(z), \cdots)$ in ${\bf T}_z$ after removing.
Thus, if
$$\{X_0={\bf 0}, X_1=u_1,X_2=u_1, \cdots\}\in \{h_m \geq Cm, \rho_0\leq N\} ,$$
then
$$\{X_{i+1}^i=z, X_{i+2}^i(z)=u_1(z), \cdots\}\in \{h_m ^i(z)\geq Cm, \rho_0^i\leq N\} .$$
On the other hand, given  a fixed paths $\{{\bf 0}, u_1,\cdots, u_j, \cdots\}$, it follows from the definition of 
 $\{z, u_{1}(z), \cdots, u_{j}(z), \cdots\}$ that
$$w\left([u_j ,u_{j+1}], k\right)= w\left([u_{j}(z), u_{j+1}(z)], i+1+k\right)\eqno{(3.24)}$$
for any positive integers $j$ and $k$.
We may focus on a finite part $\{{\bf 0}, u_1, \cdots u_l\}$ from $\{{\bf 0}, u_1, \cdots\}$.
Now if we can show that for all large $l$, and for each path $\{{\bf 0}, u_1, u_2, \cdots, u_l\}$, 
\begin{eqnarray*}
&&{\bf P} (X_1=u_1, X_2=u_2, \cdots, X_l=u_l)\\
&\leq &2^N {\bf P}\left(X_{i+2}^i(z)=u_1(z), X_{i+3}^i(z)=u_2(z),\cdots , X_{i+2+l}^i(z)=u_l(z)\,\,|\,\, w(e_z, i+1)=2\right),
\hskip 0.2cm {(3.25)}
\end{eqnarray*}
 then (3.23) will be followed by the summation  of  all possible paths
 $\{{\bf 0}, u_1, u_2, \cdots u_l\}$  for both sides  in (3.25) and 
 by letting $l\rightarrow \infty$.
 Therefore, to show (3.23),  we need  to show (3.25). 
 
 Note that
 $${\bf P} (X_1=u_1, X_2=u_2, \cdots, X_l=u_l)=\prod_{j=1}^l {\bf P}(X_j=u_j\,\,|\,\, X_{j-1}=u_{j-1},\cdots, 
 X_1=u_1)\eqno{(3.26)}$$
 and
 \begin{eqnarray*}
 &&{\bf P} (X_{i+2}^i=u_1(z), X_{i+3}^i(z)=u_2(z), \cdots, X_{i+2+l}^i(z)=u_l(z))\\
 &=&\!\!\!\!\prod_{j=1}^l {\bf P}(X_{i+1+j}^i(z)=u_j(z)\,\,|\,\, X_{i+j}^i(z)=u_{j-1}(z),\cdots ,
 X_{i+2}^i(z)=u_1(z), w(e_z, i+1)=2). \hskip 0.1cm (3.27)
 \end{eqnarray*}
  If $u_{j-1}={\bf 0}$, then
 $${\bf P}(X_j=u_j\,\,|\,\, X_{j-1}=u_{j-1},\cdots ,X_1=u_1)={w([u_{j-1}, u_{j}], j-1) \over \sum_{e} w(e, j)},
 \eqno{(3.28)}$$
 where the sum in (3.28) takes over all possible edges adjacent to the root with vertices in ${\bf T}$.
 On the other hand, if $u_{j-1}={\bf 0}$, we know that $u_{j-1}(z)=z$, then  by (3.22),
 \begin{eqnarray*}
 &&{\bf P}(X_{i+1+j}^i(z)=u_j(z)\,\,|\,\, X_{i+j}^i(z)=u_{j-1}(z),\cdots ,
 X_{i+2}^i(z)=u_1(z),w(e_z, i+1)=2)\\
 &=& {w([u_{j-1}(z), u_{j}(z)], i+j) \over \sum_{e} w(e,i+ j)+w(e_z,i+ j)}={w([u_{j-1}(z), u_{j}(z)], i+j) \over \sum_{e} w(e, i+j)+2} ,
 \hskip 4cm (3.29)
 \end{eqnarray*}
 where the sum in (3.29) takes all edges adjacent to $z$ with vertices in ${\bf T}_z$ (not including $e_z$).
 We check the numerators in the right sides of (3.28) and (3.29). If $X_1, \cdots X_{j-1}$ never visit $u_j$, then both $w([u_{j-1}, u_j],j-1]=1$ and $w([u_{j-1}(z), u_{j}(z)], i+j) =1$. Otherwise,
 by (3.24) the two numerators  are also the same. Similarly,  the two sums in the denominators 
 in the right sides of (3.28) and (3.29) are the same. 
 Therefore, if $u_{j-1}={\bf 0}$, note that $ \sum_{e} w(e, j)\geq 2$ for all $j$, so
 \begin{eqnarray*}
 && 2 {\bf P}(X_{i+1+j}^i(z)=u_j(z)\,\,|\,\, X_{i+j}^i(z)=u_{j-1}(z),\cdots ,
 X_{i+2}^i(z)=u_1(z),w(e_z, i+1)=2) \\
 &\geq & {\bf P}(X_j=u_j\,\,|\,\, X_{j-1}=u_{j-1},\cdots ,
 X_1=u_1).\hskip 6cm (3.30)
 \end{eqnarray*}
 If $u_{j-1}\neq {\bf 0}$, we do not need to consider the extra term $w(e_z,i+ j)$ in the denominator
 of the second right side of (3.29). So by the same argument of (3.30), if $u_{j-1}\neq {\bf 0}$,
 \begin{eqnarray*}
 && {\bf P}(X_{i+1+j}^i(z)=u_j(z)\,\,|\,\, X_{i+j}^i(z)=u_{j-1}(z),\cdots ,
 X_{i+2}^i(z)=u_1(z),w(e_z, i+1)=2) \\
 &=& {\bf P}(X_j=u_j\,\,|\,\, X_{j-1}=u_{j-1},\cdots ,
 X_1=u_1)\hskip 6cm (3.31)
 \end{eqnarray*}
Since we restrict $\rho_0\leq N$ and $\rho^i_0\leq N$,   walks $\{X_1,X_2,\cdots \}$
 visit the root no more than $N$ times. On the other hand, walks $\{X^i_{i+2}(z), X_{i+3}^i(z), \cdots \}$ also
visit $z$ no more than $N$ times.    This indicates that there are at most $N$ vertices $u_j$s with 
$u_j={\bf 0}$ for $1\leq j\leq l$  for the above path
$\{{\bf 0}, u_1, \cdots, u_l \}$.  Thus, (3.25) follows from (3.26)-(3.31). So does (3.23).

With (3.23), we will show Lemma 3.3.
Note that $\{h_{m}^i(z)\geq m C, \rho_0^i\leq N\}$ only depends on the weight configurations of the edges with vertices inside ${\bf T}_z$, and weight $w(e_z, i+1)$,  and the time interval $[i+2, \infty)$. 
In contrast,
on ${\cal G}(x, z, i,N, nC)$, the last move of walks $\{X_1, \cdots, X_{i}, X_{i+1}\}$ is from $x$ to $z$, but the other moves use the edges with the vertices inside $\{y: h(y)\leq h(z)-1\}$. So by (3.23), 
\begin{eqnarray*}
&&{\bf P} (h_m\geq Cm, \rho_0 \leq N)\\
&\leq & 2^N {\bf P}\left(h_m^i(z) \geq Cm, \rho_0^i\leq N\,\,|\,\, w(e_z, i+1)=2\right)\\
&\leq & 2^N{\bf P}\left(h_m^i (z) \geq Cm, \rho_0^i\leq N\,\,|\,\, {\cal G}(x,z,i,N, nC)\right).\hskip 4.8cm  {(3.32)}
\end{eqnarray*}
By  (3.21) and (3.32),
\begin{eqnarray*}
&&{\bf P} ( h_n\geq n C,\rho_0\leq N){\bf P}(h_m\geq mC,\rho_0\leq N)\\
&\leq &\sum_{i\leq n} \sum_x 2^{N}(b+1){\bf P} \left({\cal G}(x, z,i, N,nC), h_m^i(z)\geq mC, \rho^i_0\leq N\right).\hskip 3cm {(3.33)}
\end{eqnarray*}
If $i\leq n$, then
$$h_{m}^i (z)\leq h_{m+n-i}^i(z).\eqno{(3.34)}$$
By (3.33) and (3.34),
\begin{eqnarray*}
&&{\bf P} ( h_n\geq n C, \rho_0\leq N){\bf P}(h_m\geq mC,\rho_0\leq N)\\
&\leq & \sum_{i\leq n}\sum_x 2^N(b+1){\bf P} \left( {\cal G}(x,z, i, N,nC), h_m^i(z)\geq mC, \rho^i_0\leq N\right)\\
&\leq & \sum_{i\leq n} 2^N (b+1){\bf P} \left(\bigcup_x \left\{ {\cal G}(x, z,i,N, nC),h_{m+n-i}^i(z)\geq m C\right\}\right).\hskip 3.5cm (3.35)
\end{eqnarray*}
Note that for each $x$ and $i$,
$$ \left\{ {\cal G}(x,z,i, N,nC),h_{m+n-i}^i(z)\geq m C\right\}$$
implies that the walks first move to $x$ at time $i$ with $h(x) \geq nC$ and the number of  walks
$\{X_1, \cdots, X_{i}\}$  back to   the root is not more than $N$.
After that, the walks continue to move from $x$ to $z$. After this move, the walks move
 inside subtree ${\bf T}_z$.  So $i$ is a cut time and $X_{i}$ is a cut vertex with $h(X_i) \geq nC$.
 Therefore, together with $h_{n+m-i}^i(z) \geq mC$,
$ \left\{ {\cal G}(x, z, i,N,  nC),h_{m+n-i}^i(z)\geq m C\right\} $  implies that $\{h_{n+m+1} \geq (n+m)C+1, \rho_0\leq N\}$ occurs. In other words,
$$ \left\{ {\cal G}(x, z, i,N,  nC),h_{m+n-i}^i(z)\geq m C\right\}\subset
\{h_{n+m+1}\geq (n+m)C+1, \rho_0\leq N\}.\eqno{(3.36)}$$
Therefore,
$$\bigcup_{x} \left\{ {\cal G}(x,z, i, N,nC),h_{m+n-i}^i(z)\geq m C\right\}\subset
\{h_{n+m+1}\geq (n+m)C+1,\rho_0\leq N\}.\eqno{(3.37)}$$

Finally, by (3.35) and (3.37),
\begin{eqnarray*}
&&{\bf P} ( h_n\geq n C,\rho_0\leq N){\bf P} ( h_m\geq m C,\rho_0\leq N)\\
&\leq & 2^N (b+1)n {\bf P} ( h_{n+m+1}\geq (n+m)  C+1,\rho_0\leq N).\hskip 5cm {(3.38)}
\end{eqnarray*}
Therefore, Lemma 3.3 follows from (3.38). $\Box$\\

We let 
$$a_n=-\log {\bf P}(h_n \geq (T+\epsilon) n, \rho_0\leq N).\eqno{(3.39)}$$
We may take $\epsilon$ small such that $T+\epsilon< 1$.  By Lemma 3.3, for any $n$ and $m$
$$a_{n+m+1}\leq a_n+ a_m + \log n + N \log 2+\log(b+1).\eqno{(3.40)}$$
By (3.40) and a standard subadditive argument (see (II.6) in Grimmett (1999)), we have the following lemma.\\

{\bf Lemma 3.4.} {\em For $c=2$ and any $N>0$ and $b \geq 2$, there exists $0\leq  \alpha(N) <\infty$ such that}
 $$\lim_{n\rightarrow \infty} {-1\over n} \log {\bf P}(h_n  \geq (T+\epsilon)n, \rho_0\leq N) =\inf_{n} \left\{ {-1\over n} \log {\bf P}(h_n  \geq (T+\epsilon)n, \rho_0\leq N)\right\}=\alpha(N). $$

It follows from the definition and Lemma 3.4 that $\alpha(N)$ is a non-negative decreasing  sequence in $N$. 
Thus,  there exists a finite number $\alpha\geq 0$ such that
$$\lim_{N\rightarrow \infty} \alpha(N)=\alpha.\eqno{(3.41)}$$
By (3.41) and Lemma 3.4,  for each $N$,
$$\alpha \leq \alpha(N)\leq {-1\over n} \log {\bf P}(h_n  \geq (T+\epsilon)n, \rho_0\leq N).\eqno{(3.42)}$$
On the other hand, note that the walk is transient, so $\rho_0 < \infty$. Thus, for any fixed $n$, 
$$\lim_{N\rightarrow \infty }{-1\over n} \log {\bf P}(h_n  \geq (T+\epsilon)n, \rho_0\leq N)={-1\over n} \log {\bf P}(h_n  \geq (T+\epsilon)n)\eqno{(3.43)}$$
By (3.42) and (3.43),
$$\alpha\leq \liminf_n{-1\over n} \log {\bf P}(h_n  \geq (T+\epsilon)n)\eqno{(3.44)}$$
Note that for each $N$,
$$\limsup_{n} {-1\over n} \log {\bf P}(h_n  \geq (T+\epsilon)n)\leq \lim_{n\rightarrow \infty} {-1\over n} \log {\bf P}(h_n  \geq (T+\epsilon)n, \rho_0\leq N)=\alpha(N).$$
So for each $\delta >0$ we may take $N$ large such that

$$ \limsup_n{-1\over n} \log {\bf P}(h_n  \geq (T+\epsilon)n)\leq \alpha(N) \leq \alpha+\delta.\eqno{(3.45)}$$
We summarize (3.44) and (3.45) as the following lemma.\\

{\bf Lemma 3.5}.  {\em For $c=2$ and any $b \geq 2$, there exists a constant $\alpha\geq 0$ such that}
$$ \lim_{n\rightarrow \infty}{-1\over n} \log {\bf P}(h_n  \geq (T+\epsilon)n)= \alpha.$$


\section{ Proof of Theorem 1.}
Note that for $\epsilon < 1- T$, and for all large $n$,
$$\left({b\over b+1}\right)^n \leq {\bf P}(h(X_{i+1})> h(X_i)\mbox{ for } 0\leq i\leq n) \leq {\bf P}(h(X_n)\geq n(T+\epsilon)).\eqno{(4.1)}$$
By (4.1),
$$\limsup_{n\rightarrow \infty} {-1\over n} \log {\bf P} ( h(X_n)\geq n(T+\epsilon))< \infty.\eqno{(4.2)}$$
 Note also that
\begin{eqnarray*}
&&{\bf P} ( h(X_n)\geq n(T+\epsilon))\\
&\leq&{\bf P} ( h(X_n)\geq n(T+\epsilon), h_n(N,M) \geq n(T+\epsilon/2))+{\bf P} ( h(X_n)- h_n(N,M) \geq   n\epsilon/2).\hskip 0.3cm {(4.3)}
\end{eqnarray*}
By Lemmas 2.3 and 2.4, for $\epsilon >0$, we select $N$ and $M$ such that
$${\bf P} ( h(X_n)\geq n(T+\epsilon))\leq {\bf P} ( h_n(N,M) \geq n(T+\epsilon/2))+C_1\exp(-C_2 n).\eqno{(4.4)}$$
For $N$ and $M$ in (4.4), we may require that $N \geq N_0$ and $M \geq M_0$ for $N_0$ and $M_0$
in Lemma 3.2.
By (4.4) and Lemma 3.2, there exist $C_i=C_i(\epsilon, N,M)$  for $i=3,4$ such that
$${\bf P} ( h(X_n)\geq n(T+\epsilon))\leq C_3\exp(-C_4 n).\eqno{(4.5)}$$
By (4.5), for $\epsilon >0$, 
$$0<  \liminf_{n\rightarrow \infty} {-1\over n} \log {\bf P} ( h(X_n)\geq n(T+\epsilon)).\eqno{(4.6)}$$
 
It remains for us to show the existence of the limit in Theorem 1.
We use a similar proof in  Lemma 3.3 to show it.
Let ${\cal F}(x, k, n)$ be the event that 
 $h(X_i) < n(T+\epsilon)$ for $ i=1,\cdots,k-1$,  $h(X_{k})\geq n(T+\epsilon)$ and $h(X_k)=x$ for $k\leq n$.
 Thus,
$${\bf P}( h(X_n) \geq n(T+\epsilon))\leq  \sum_{k\leq n} \sum_{x\in {\bf T}} {\bf P} ({\cal F}(x, k, n))\eqno{(4.7)} $$
Note that ${\cal F}(x, k,  n)$ depends on  finite step walks $\{X_0, \cdots, X_{k}\}$. We need to couple
the remaining walks $\{X_{k+1}, X_{k+2}, \cdots\}$ such that $k$ is a cut time.
 Let  ${\cal Q}(x,k)$   be the event that $X_{k}=x$  and $\{X_t\}$ will stay inside $ {\bf T}_x$ but never returns
to $x$ for $t > k$.  Since the walks are transient,  we may let
$${\bf P} ({\cal Q}({\bf 0}, 0))=\nu >0.\eqno{(4.8)}$$
Let  $e_x$ denote the edge with vertices $x$ and $w$ for $h(w) =h(x)-1$. We know
 that ${\cal Q}(x, k)$ depends on  initial weight $w(e_x, k)$, and   the weights in the  edges  with the vertices in ${\bf T}_x$,
respectively.
Therefore, by the same discussion of (3.23) in Lemma 3.3,
$$2{\bf P} ({\cal Q}(x,k)\,\,\, |\,\,\,{\cal F}(x, k, n))\geq \left({b+2\over b}\right){\bf P} ({\cal Q}(x,k)\,\,\, |\,\,\,{\cal F}(x, k, n)) \geq \nu .\eqno{(4.9)}$$
Thus, by (4.7) and (4.9),
\begin{eqnarray*}
&&{\bf P} (h(X_n)\geq n(T+\epsilon))\\
&\leq &\sum_{x\in {\bf T}}\sum_{k\leq n} {\bf P} \left(  {\cal F}(x, k,n) \right){\bf P} ({\cal Q}(x,k)\,\,\, |\,\,\,{\cal F}(x,k,n))\left({b+2\over b}\right)\nu^{-1}\\
&\leq &  2\nu^{-1}\sum_{k\leq n}{\bf P} \left(  \bigcup_{x\in {\bf T} }{\cal F}(x,k,n)\cap  {\cal Q}(x,k)\right).\hskip 7cm {(4.10)}
\end{eqnarray*}
If
$ {\cal F}(x,k,n) \cap  {\cal Q}(x,k)$
occurs, it implies that  the walks move to $x$ at $k \leq n$ with $h(x) \geq n(T+\epsilon)$. 
After that, the walks continue to move
inside ${\bf T}_x$ from $x$ and never return to $x$. This implies that 
$k$ is a cut time and $X_{k}$ is a cut vertex with $h(X_{k}) \geq n(T+\epsilon)$. 
So for $0\leq k\leq n$ and for each $x$,
$${\cal F}(x,k,n) \cap {\cal Q}(x,k)\subseteq \{h_{k}\geq n(T+\epsilon)\}. \eqno{(4.11)}$$
Thus,
$$\bigcup_{x\in {\bf T}}{\cal F}(x,k,n) \cap {\cal Q}(x,k)\subseteq \{h_{k}\geq n(T+\epsilon)\}. \eqno{(4.12)}$$
Note that for $0\leq k\leq n$,
$$h_{k} \leq h_n.\eqno{(4.14)}$$
By (4.10)-(4.14),
$${\bf P} (h(X_n)\geq n(T+\epsilon))\leq 2\nu^{-1} n {\bf P}(h_{n}\geq n(T+\epsilon)).\eqno{(4.15)}$$

On the other hand, we suppose that $h_{n}\geq n(T+\epsilon)$. 
Note that if $\tau_k\leq n\leq \tau_{k+1}$, then by (3.18),
$$h_n=h\left(X_{\tau_k}\right)\leq h(X_n).\eqno{(4.16)} $$
By (4.16),
$${\bf P} (h_{n}\geq n(T+\epsilon))\leq  {\bf P} (h(X_n) \geq n(T+\epsilon)).\eqno{(4.17)}$$
Now we are ready to show Theorem 1. \\

{\bf Proof of Theorem 1.} Together with (4.15), (4.17), and Lemma 3.5, 
$$\lim_{n\rightarrow \infty} {1\over n} \log {\bf P} ( h(X_n)\geq n(T+\epsilon))=\alpha.\eqno{(4.18)}$$
By (4.2) and (4.6),
$$0< \alpha < \infty.\eqno{(4.19)}$$
Therefore, Theorem 1 follows from (4.18) and (4.19). $\Box$

\section{Proof of Theorem 2.} Similarly, we define the same cut times
$\tau_i$ that we defined for the linearly reinforced random walk.
We have $\left(\tau_{k+1}-\tau_k, h(Y_{\tau_{k+1}})-h(Y_{\tau_k})\right)$ as an i.i.d. sequence.
We can also follow  
Durrett, Kesten, and Limic's (2002)  Lemmas 7 and 8 to show that there exist $C_i$ for $i=1,2$ such that,
for each $k\geq 1$,
$${\bf P}( \tau_{k+1}-\tau_k\geq m)\leq C_1\exp(-C_2m)\eqno{(5.1)}$$
and 
$${\bf P}( h(Y_{\tau_{k+1}})-h(Y_{\tau_k})\geq m)\leq C_1\exp(-C_2m).\eqno{(5.2)}$$
By (5.1) and (5.2),
similar to our approach the linearly reinforced random walk, we set 
$$S_n= \sum_{k=1}^n (\tau_k-\tau_{k-1})\mbox{ and } K_n= \sum_{k=1}^n \left(h(Y_{\tau_{k}})-h(Y_{\tau_{k-1}})\right).\eqno{(5.3)}$$
By the law of large numbers, 
$$\lim_{n\rightarrow \infty}{S_n \over n}=A\mbox{ and }\lim_{n\rightarrow \infty}{K_n \over n}=B.\eqno{(5.4)}$$
With these observations, Theorem 2  can follow from the exact proof of Theorem 1. In fact,
we may not need to truncate $\tau_i$ to $\tau_i(N,M)$ as we did for Theorem 1, since we can use (5.1) and (5.2)
directly. $\Box$\\

\section { Proof of Theorem 3.}

Now we need to estimate ${\bf P}( h(Y_n)\leq n(S-\epsilon))$. Let 
$$L_n=\max \{i, \tau_i \leq n\}$$
and let
$$h_n=\sum_{i=1}^{L_n} \left(h\left(Y_{\tau_{i}}\right)-h\left(Y_{\tau_{i-1}}\right)\right) \mbox{ if }L_n\geq 1\mbox{ and }  h_n=0 \mbox{ if } L_n=0.\eqno{(6.1)}$$

By (1.7), (5.3), and an argument similar to (3.6), we have
$$ {B\over A}=S.\eqno{(6.2)}$$
Since $h_n \leq h(Y_n)$, by (5.1)
\begin{eqnarray*}
&&{\bf P}( h(Y_n)\leq n(S-\epsilon B))\\
&\leq & {\bf P}(  h_n \leq n(S-\epsilon B))\\
&\leq &  {\bf P}\left( \sum_{i=1}^{L_n}  \left(h\left(Y_{\tau_{i}}\right)-h\left(Y_{\tau_{i-1}}\right)\right) \leq n(S-\epsilon B\right)+{\bf P}( \tau_1 > n)\\
&\leq &  {\bf P}\left( \sum_{i=1}^{L_n}  \left(h\left(Y_{\tau_{i}}\right)-h\left(Y_{\tau_{i-1}}\right)\right) \leq n(S-\epsilon B\right)+C_1\exp(-C_2 n).\hskip 3cm {(6.3)}
\end{eqnarray*}
We split
\begin{eqnarray*}
&&{\bf P}\left( \sum_{i=1}^{L_n}  (h(Y_{\tau_{i}})-h(Y_{\tau_{i-1}})) \leq n(S-\epsilon B)\right)\\
&\leq &{\bf P}\left( \sum_{i=1}^{L_n}  (h(Y_{\tau_{i}})-h(Y_{\tau_{i-1}})) \leq n(S-\epsilon B), L_n \geq n(SB^{-1} 
-\epsilon /2)\right)\\
&&+{\bf P}\left(   L_n < n(SB^{-1} -\epsilon/2)\right)\\
&=&I+II.
\end{eqnarray*}
We estimate $I$ and $II$ separately:
\begin{eqnarray*}
I&=&{\bf P}\left( \sum_{i=1}^{L_n}  (h(Y_{\tau_{i}})-h(Y_{\tau_{i-1}}) )\leq n(S-\epsilon B), L_n \geq n(SB^{-1} -\epsilon/2)\right)\\
&\leq &{\bf P}\left( \sum_{i=1}^{n(SB^{-1} -\epsilon/2)}  (h(Y_{\tau_{i}})-h(Y_{\tau_{i-1}})) \leq n(S-\epsilon B)\right).\hskip 4.5cm (6.4)
\end{eqnarray*}
Note that
$${\bf E}\left(\sum_{i=1}^{n(SB^{-1} -\epsilon/2)}  (h(Y_{\tau_{i}})-h(Y_{\tau_{i-1}}))\right)= n(S -\epsilon B/2).\eqno{(6.5)}$$
Note also that by (5.2), 
$\{h(Y_{\tau_{i}})-h(Y_{\tau_{i-1}})\}$ is an i.i.d. sequence with an exponential tail for $k\geq 2$, so by Lemma 2.2
there exist $C_i=C_i(\epsilon, B)$ for $i=3,4$ such that
$$I\leq C_3\exp(-C_4 n).\eqno{(6.6)}$$
Also, by (6.2),
$$II={\bf P}\left(   L_n < n(SB^{-1} -\epsilon/2)\right)={\bf P} \left(\sum_{i=1}^{n(SB^{-1} -\epsilon/2)} 
(\tau_i-\tau_{i-1})\geq n\right) ={\bf P} \left(\sum_{i=1}^{n(A^{-1} -\epsilon/2)} 
(\tau_i-\tau_{i-1})\geq n\right).$$
Note that
$${\bf E}\sum_{i=1}^{n(A^{-1} -\epsilon/2)} (\tau_i-\tau_{i-1})= n(1-\epsilon A/2).\eqno{(6.7)}$$
Note also that by (5.1), $\{\tau_i-\tau_{i-1}\}$ is an i.i.d. sequence with an exponential tail for $k\geq 2$, 
so by Lemma 2.2,
there exist $C_i=C_i(\epsilon, B)$ for $i=5,6$ such that
$$II \leq C_5 \exp(-C_6 n).\eqno{(6.8)}$$
Together with (6.3), (6.4), (6.6),  and (6.8), there exist $C_i=C_i(c, \epsilon, B)$ for $i=7,8$ such that
$${\bf P}( h(Y_n)\leq n(S-\epsilon))\leq  C_7 \exp(-C_8 n).\eqno{(6.9)}$$
From (6.9), 
$$0< \liminf {-1\over n} \log {\bf P}( h(Y_n)\leq n(S-\epsilon)).\eqno{(6.10)}$$
If the walks repeatedly move in the edge connecting the origin in $n$ times, we have the probability
$C^n$ for a positive constant $C=C(b)$. Thus, for $\epsilon < S$ and for all large $n$,
$$ C^n\leq {\bf P}(h(Y_n)\leq 1) \leq {\bf P}(h(Y_n)\leq n(S-\epsilon)).\eqno{(6.11)}$$
So for $\epsilon < S$,
$$\limsup {-1\over n} \log {\bf P}( h(Y_n)\leq n(S-\epsilon))< \infty.\eqno{(6.12)}$$
Therefore, Theorem 3 follows from (6.10) and (6.12). $\Box$\\

{\bf Acknowledgments.} The author would like to thank a  careful referee for his many valuable comments, which corrected a few mistakes and improved the quality of this paper. He would also
like to thank Takei M. for pointing out  a few typos and a few references.\\

\newpage
\begin{center}
{\bf References}
\end{center}
1. Aidekon, E. (2008). Transient random walks in random environment on a Galton Watson tree.
{\em Probab. Theory Related Fields} {\bf 142}, 525--559.\\
2. Collevecchio, A. (2006). Limit theorems for reinforced random walks on certain trees.
{\em Probab. Theory Relate  Fields} {\bf 136}, 81--101.\\
3. Davis, B. (1990). Reinforced random walks. {\em Probab. Theory Relate  Fields}.  {\bf 84}, 203-229.\\
4. Diaconis, P. (1988). Recent progress on de Fietti's notions of exchangeability. {\em Bayesian statistics} {\bf 3},
115--125, Oxford Univ. Press.\\
5. Engel, E. (1997). {\em Sperner Theory}. Cambridge Univ. Press, New York.\\
6. Durrett, R., Kesten, H., and Limic, V. (2002). Once edge-reinforced random walk. {\em Probab. Theory Related Fields}
{\bf 122}, 567--592.\\
7. Grimmett, G. (1999). {\em Percolation}. Springer-Verlag,  New York.\\
8. Pemantle, R. (1988). Phase transition in reinforced random walks and RWRE on trees. {\em Ann. Probab}. {\bf 16}, 
1229--1241.\\
9. Sabot, C and Tarres, P (2012).  Edge-reinforced random walk, vertex-reinforced jump process and the supersymmetric hyperbolic sigma model.  {\em arXiv} 1111.3991v3.\\
10. Takeshima, M. (2000) Behavior of 1-dimensional reinforced random walk. {\em  Osaka J. Math}. {\bf 7}, 355-372. \\

\end{document}